\documentclass[12pt]{article}
\usepackage{latexsym}

\newcommand{\duk}{\noindent {\bf Proof. }}
\newcommand{\kduk}{\hfill $\Box$\bigskip}
\newcommand{\R}{\mathbf{R}}
\newcommand{\N}{\mathbf{N}}
\newcommand{\Z}{\mathbf{Z}}
\newcommand{\Q}{\mathbf{Q}}

\title{Non-holonomicity of the sequence $\log 1, \log 2, \log 3, \dots$}
\author{Martin Klazar\thanks{Institute for Theoretical Computer Science and Department of Applied Mathematics, 
Faculty of Mathematics and Physics, Charles University, Malostransk\'e n\'am\v est\'\i\ 25, 118 00 Praha, Czech Republic. 
ITI is supported by the project 1M0021620808 of the Czech Ministry of Education.
Email: {\tt klazar@kam.mff.cuni.cz}}}
\date{\today}

\begin{document}

\maketitle
\begin{abstract}
Gerhold conjectured and proved conditionally that $(\log n)_{n\ge 1}$ is not a holonomic sequence. Flajolet, Gerhold and Salvy
gave a proof using an analytic machinery. We give a simple proof.
\end{abstract}

A sequence of complex numbers $(a_n)_{n\ge 1}\subset\mathbf{C}$ is called {\em holonomic} (P-recursive, D-finite) if there 
exist polynomials $p_i(x)\in\mathbf{C}[x]$, $i=0,1,\dots,d$, which are not all zero and
\begin{equation}\label{Prec}
p_d(n)a_{n+d}+\cdots+p_1(n)a_{n+1}+p_0(n)a_n=0\ \mbox{ for every integer $n\ge 1$.}
\end{equation}
It can be proved that if $(a_n)_{n\ge 1}\subset\R$, the polynomials $p_i$ can be taken with
real coefficients as well, and the same is true for the field $\Q$ of rationals (which is the case in enumerative applications). 
Therefore, since we consider the real sequence $(\log n)_{n\ge 1}$, without loss of generality we can restrict the polynomials 
in (\ref{Prec}) to $\R[x]$ (but not to $\Q[x]$). For the importance of holonomic
sequences in combinatorial enumeration and for further information see the book \cite{stan} of Stanley. 

Gerhold \cite{gero} used Schanuel's conjecture to prove that $(\log n)_{n\ge 1}$
is not holonomic and asked about an unconditional proof. Such proof was recently given by Flajolet, Gerhold and Salvy 
in \cite{flaj_al} and is based on the equivalence saying that $(a_n)_{n\ge 1}$ is holonomic iff the generating function 
$\sum_{n\ge 1}a_nx^n$ satisfies a linear differential equation with polynomial coefficients. The authors further apply a
theorem characterizing asymptotic behavior of solutions of such equations near singularity and an 
Abelian theorem on transfer of the asymptotics of a sequence to its generating function. 

In this note we give a simple proof of their result that $(\log n)_{n\ge 1}$ is not holonomic, that is, the sequence
$a_n=\log n$ satisfies no nontrivial relation (\ref{Prec}). We prove a somewhat stronger statement. (Recall that all polynomials 
and rational functions here have real coefficients.)

\bigskip\noindent
{\bf Theorem. }{\em 
Let $b_0=0<b_1<\dots<b_d$ be real numbers and $p_0(x),p_1(x),\dots,p_d(x)$ be polynomials, not all of them zero. 
The function
$$
F(x)=\sum_{i=0}^dp_i(x)\log(x+b_i)
$$ 
has only finitely many positive real zeros. In particular, $(\log n)_{n\ge 1}$ is not holonomic.
}

\smallskip\noindent
\duk
First we suppose  that $F(x)=0$ for every $x>0$ and derive a 
contradiction. We may assume that $p_0$ is nonzero. Let $m\ge 0$ be minimum such that $p_0^{(m)}(0)=c\ne 0$.
Differentiating $m$ times the equation $F(x)=0$, we obtain that
$$
F^{(m)}(x)=u(x)+s(x)=\bigg(\sum_{i=1}^dp_i^{(m)}(x)\log(x+b_i)+p_0^{(m)}(x)\log x\bigg)+s(x)=0
$$
for every $x>0$, where in the rational function $s(x)$ we have collected all terms $q(x)/(x+b_i)^k$ with $q(x)\in\R[x]$ arising 
in the differentiation ($s(x)=0$ if there is no such term). Note that $0$ is not a pole of $s(x)$ (even if it were, the following 
argument would still work) and that $s(x)$ has no real positive pole. $F^{(m)}(x)$ is 
defined for every $x>0$. For $x\to 0^+$ we have $|u(x)|\sim |c|\cdot\log(1/x)\to\infty$ but $s(x)$ is bounded.  This contradicts 
the equality $u(x)+s(x)=0$ for all $x>0$.

Now let $F(x_j)=0$ for an infinite sequence $0<x_1<x_2<\dots$ (or $x_1>x_2>\dots >0$). We consider the derivatives 
$F',F'',\dots,F^{(M+1)}$ where $M$ is the largest degree of $p_i$. Clearly, $F^{(M+1)}$ is a rational function. Note 
that all these derivatives are defined for every $x>0$.  Applying Rolle's theorem (it says that if $f$ 
is continuous on $[a,b]$, has derivative on $(a,b)$, and $f(a)=f(b)$, then $f'(c)=0$ for some $c\in(a,b)$), we obtain 
that each derivative of $F$ has infinitely many positive real zeros as well. Since every nonzero rational function has only 
finitely many zeros, $F^{(M+1)}$ must be identically zero. Thus $F^{(M)}$ is constant and identically zero. 
The same holds for all previous derivatives and we conclude that $F$ is identically zero (for $x>0$). But we have shown that this 
is not possible. 
\kduk

It can be shown very quickly that (\ref{Prec}) with $a_n=\log n$ cannot be satisfied 
with {\em rational} polynomials $p_i$. Let us suppose there is such relation. We multiply it by an integer to clear out 
denominators and get integral polynomials. Then we apply exponential function and get the relation
\begin{equation}\label{fraction}
\frac{\prod_{i\in I}(n+i)^{q_i(n)}}{\prod_{i\in J}(n+j)^{q_i(n)}}=1
\end{equation}
where $I$ and $J$ are disjoint subsets of $\{0,1,\dots,d\}$, not both empty (for empty
set the corresponding product is defined to be $1$), and every $q_i(x)$ is a nonzero integral polynomial with 
positive leading coefficient. We take the largest element of $I\cup J$, let it be $a\in I$. For sufficiently large $n\in\N$ all 
values $q_i(n)$ are positive integers. We take such an $n$ with the property that $n+a=p$ is a prime. 
Then the numerator 
and the denominator in (\ref{fraction}) are positive integers such that the former is divisible by $p$ but the latter is not. They 
cannot cancel out as required by the $1$ on the right side and we have a contradiction.

\bigskip\noindent
{\bf Concluding remarks.} Almost all proofs of non-holonomicity of sequences work with their generating functions, the exceptions 
being \cite{gero} and the present note. Our proof suggests another analytic method for non-holonomicity
proofs that is applicable when $a_n=f(n)$ for sufficiently smooth function $f(x)$. Instead of the generating function, 
one can take directly (\ref{Prec}) which implies that 
$$
F(x)=\sum_{i=0}^d p_i(x)f(x+i)
$$
has infinitely many positive  zeros. By Rolle's  theorem, so has each derivative of $F$. Then one can try to bring this 
to contradiction. For example, this approach should work for $a_n=n^{\alpha}$, $\alpha\in\R\backslash\Z$ (another sequence 
considered in \cite{flaj_al} and \cite{gero}).

\end{document}